\documentclass[a4paper, 11pt]{article}
\usepackage{comment}
\usepackage[
margin = 1.5cm,
includefoot,
footskip = 30pt,
]{geometry}\usepackage{enumitem}
\usepackage[english]{babel}
\usepackage[utf8]{inputenc}
\usepackage{amsmath}
\usepackage{amsfonts}\usepackage{ amssymb }\usepackage{graphicx}\usepackage{float}
\usepackage{amsthm}\theoremstyle{plain}
\newtheorem{theorem}{Theorem}
\newtheorem{proposition}{Proposition}
\newtheorem{conjecture}{Conjecture}
\newtheorem{lemma}{Lemma}
\newtheorem{remark}{Remark}

\usepackage{mathtools}
\usepackage{babel}

\title{On the class-breadth conjecture\footnote{
This research was supported by the Russian Science Foundation under grant no. 22-11-00075.}}
\author{Alexander Skutin
\\{\small Faculty of Mechanics and Mathematics of Lomonosov Moscow State University}\\{\small Moscow Center for Fundamental and Applied Mathematics}\\{\small a.skutin@mail.ru}}
\date{}

\begin{document}
\maketitle
\begin{abstract}
    The class-breadth conjecture of Leedham-Green, Neumann and Wiegold states that for each $p$-group, $cl(G)\leq b(G) + 1$, where $cl(G)$, $b(G)$ denote the nilpotency class and the breadth of $G$. While several counter-examples to this conjecture have been found for $p = 2$, it is still open in general for $p > 2$. This article is dedicated to the general case $p > 2$ of the conjecture.
\end{abstract}

\section{Introduction}

Let $G$ be a finite $p$-group. The breadth $b(x)$ of an element $x$ of a group $G$ is defined as $b(x) := \log_p|G:Z_G(x)|$, where $Z_G(x)$ denotes the centralizer of $x$ in $G$. The breadth $b(G)$ of $G$ is defined as $\displaystyle b(G) := \max_{g\in G}b(g)$.

C. Leedham-Green, P. Neumann and J. Wiegold in \cite{1} stated the following conjecture: $cl(G)\leq b(G) + 1$, for all $p$-groups $G$ ($cl(G)$ denotes the nilpotency class of $G$).

W. Felsch, J. Neubüser and W. Plesken in \cite{2} provide the series of counter-examples to the conjecture for $p = 2$. Further counter-examples were found by B. Eick, M. Newman and E. O'Brien in \cite{3}, in \cite{3} they also proved that the analogue of the conjecture for finite dimensional nilpotent Lie algebras is false for all finite fields, i.e., that for every finite field $\mathbb{F}$ there exists a finite dimensional nilpotent Lie algebra $\mathfrak{g}$ over $\mathbb{F}$ such that $cl(\mathfrak{g}) > b(\mathfrak{g}) + 1$, where $\displaystyle b(\mathfrak{g}) := \max_{a\in\mathfrak{g}}(\dim \mathfrak{g} - \dim Z_{\mathfrak{g}}(a))$, $Z_{\mathfrak{g}}(a) := \{x\in\mathfrak{g}\:\vert\: [a, x] = 0\}$. The conjecture has been proved in some special cases, see, for example: if $b(G)\leq 4$ in \cite{7, 8}, if $b(G)\leq p + 1$ in \cite{4}, if $G$ is a $p$-group of maximal class in \cite{6}, if $G$ is metabelian in \cite{1}, if $G$ is not covered by its two-step centralizers in \cite{1}. See also \cite{3} for more classes of $p$-groups satisfying the conjecture. A more detailed description of the results related to the class-breadth conjecture can be found in the introduction of \cite{2}.

This article is dedicated to the general case $p > 2$ of the conjecture.

The author thanks A. M. Maksaev for the introduction to the topic by suggesting the problem \cite[Problem 4.69]{9} and Y. V. Novikova for discussions that led to this article.

The author thanks the Theoretical Physics and Mathematics Advancement Foundation “BASIS”.

\section{Statements and proofs of the main lemmas}

Consider any $p$-group $G$. The commutator $[x, y], x, y\in G$ is $x^{-1}y^{-1}xy$. For each subgroups $G_1, G_2$ of $G$, denote by $[G_1, G_2]$ the subgroup of $G$ which is generated with the set of elements $\{[g_1, g_2]\:\vert\: g_1\in G_1, g_2\in G_2\}$. For each element $g\in G$ and each subgroup $H$ of $G$, denote $|G : H| := |G|/|H|$, $b_H(g) := \log_p|H : H\cap Z_G(g)|$, $[g, H] := \{[g, h]\:\vert\: h\in H\}$. So $b(g) = b_G(g)$, $b_H(g) = \log_p|[g, H]|$. $\gamma_1(G) := G$, $\gamma_{i + 1}(G) := [\gamma_i(G), G]$, $i\geq 1$ -- denotes the lower central series of $G$, $cl(G) := \min\{n\:\vert\: \gamma_{n + 1}(G) = \{e\}\}$ -- denotes the nilpotency class of $G$. For each subset $X$ of $G$, denote by $\langle X\rangle$ the subgroup of $G$ which is generated with the set $X$. By a maximal subgroup $H$ of $G$ denote any maximal under the inclusion subgroup $H\not= G$ of $G$, note that $|G : H| = p$ and $H$ is a normal subgroup of $G$ for each maximal subgroup $H$ of $G$. By a central subgroup $H$ of $G$ denote any subgroup $H\subseteq Z(G)$.

\begin{lemma}
    Let $G$ be a $p$-group. Consider any normal subgroups $C_1, C_2$ of $G$ such that $C_1\subsetneq C_2$. Then there exists a normal subgroup $C_3$ of $G$ such that $C_1\subseteq C_3\subseteq C_2$ and $|C_2 : C_3| = p$.
\end{lemma}

\begin{proof}
    Consider the subgroup $C_1[C_2, G] := \{xy\:\vert\: x\in C_1, y\in [C_2, G]\}$ of $G$ (it's a subgroup, because $C_1, C_2, [C_2, G]$ are normal subgroups of $G$). Easy to check that $C_1[C_2, G]$ is a normal subgroup of $G$ and $C_1[C_2, G]\subseteq C_2$. From $C_1\subsetneq C_2$ it follows that $\{e\}\subsetneq C_2/C_1\subseteq G/C_1$, therefore, since $G/C_1$ is nilpotent, $C_1[C_2, G] = [C_2/C_1, G/C_1]C_1 \subsetneq C_2$. Consider any maximal under the inclusion proper subgroup $C_3$ of $C_2$ such that $C_1[C_2, G]\subseteq C_3$. Then (from the maximality of $C_3$ in $C_2$) $|C_2 : C_3| = p$. From $[C_3, G] \subseteq [C_2, G]\subseteq C_1[C_2, G]\subseteq C_3$ and $C_1\subseteq C_1[C_2, G]\subseteq C_3$ it follows that $C_3$ is a normal subgroup of $G$ that satisfies $C_1\subseteq C_3\subseteq C_2$ and $|C_2 : C_3| = p$.
\end{proof}

\begin{lemma}
    Let $G$ be a $p$-group. Consider any normal subgroups $H, C_1, C_2$ of $G$ such that $C_1\not= C_2$, $C_1, C_2\subseteq H$, and $|H : C_1| = |H : C_2| = p$. Then each subgroup $C$ of $H$ that satisfies $C_1\cap C_2\subseteq C$ is a normal subgroup of $G$.
\end{lemma}

\begin{proof}
    For each $1\leq i\leq 2$ consider the homomorphism of factorization $\pi_i : G\to G/C_i$. From $|H : C_1| = |H : C_2| = p$ it follows that $\pi_i(H)$ is a normal subgroup of $\pi_i(G)$ of order $p$, therefore $\pi_i(H)\subseteq Z(\pi_i(G))$ and $[H, G]\subseteq C_i$, $\forall 1\leq i\leq 2$. So $[C, G]\subseteq [H, G]\subseteq C_1\cap C_2\subseteq C$ and $C$ is a normal subgroup of $G$.
\end{proof}

\begin{lemma}
    Let $G$ be a $p$-group. Consider any normal subgroup $H$ of $G$ such that the group $H/[H, G]$ is non-cyclic. Then there exist at least two maximal subgroups $C_1\not= C_2$ of $H$ which are normal in $G$.
\end{lemma}

\begin{proof}
    From the fact that $H/[H, G]$ is non-cyclic it follows that there are at least two maximal subgroups $A_1\not= A_2$ of $H/[H, G]$. Let $\pi : G\to G/[H, G]$ denote the homomorphism of factorization. Then $C_1 := \pi^{-1}(A_1)\not= \pi^{-1}(A_2):= C_2$ are maximal subgroups of $H$ which are normal in $G$, because $[C_i, G]\subseteq [H, G]\subseteq C_i$, $\forall 1\leq i\leq 2$.
\end{proof}

\begin{lemma}
    Let $G$ be a $p$-group. Consider any normal subgroup $H$ of $G$ such that the group $H/[H, G]$ is non-cyclic. Let $C_1, C_2\subseteq H$ be any two maximal subgroups of $H$ that are normal in $G$. Then there exists a maximal subgroup $C$ of $H$ such that $C$ is a normal subgroup of $G$, $C\nsubseteq C_1\cup C_2$, and $C\cap C_1 = C\cap C_2$.
\end{lemma}

\begin{proof}

Consider the case $C_1\not= C_2$. Let $C$ be any subgroup of $H$ such that $C\cap C_1 = C\cap C_2 = C_1\cap C_2$ and $|H : C| = p$ (such a group exists, since $H/(C_1\cap C_2)\cong \mathbb{Z}_p\oplus\mathbb{Z}_p$). From Lemma 2 it follows that $C$ is a normal subgroup of $G$. The group $C$ is a maximal subgroup of $H$ that satisfies $C\nsubseteq C_1\cup C_2$ and $C\cap C_1 = C\cap C_2$, because: $|H : C| = p$, $C\cap C_1 = C\cap C_2 = C_1\cap C_2$, and $|H : C_1\cap C_2| = p^2$.

Consider the case $C_1 = C_2$. From Lemma 3 there exists a maximal subgroup $C\not= C_1$ of $H$ which is normal in $G$. Then $C$ satisfies the conditions of Lemma 4.
\end{proof}

\begin{lemma}
    Consider any prime number $p > 2$ and a $p$-group $G$. Consider any two proper subgroups $C_1, C_2$ of $G$. Let $b(g) = 0$, $\forall g\in G\setminus (C_1\cup C_2)$. Then $G$ is an abelian group.
\end{lemma}

\begin{proof}
    From the fact that for $p > 2$ each $p$-group cannot be covered by any of its three proper subgroups it follows that $G = \langle G\setminus (C_1\cup C_2)\rangle\subseteq \langle \{g\in G\:\vert\: b(g) = 0\}\rangle\subseteq Z(G)$, i.e. that $G$ is an abelian group.
\end{proof}

\section{Strengthened class-breadth conjecture and its proof for a special class of $p$-groups}

\begin{conjecture} Consider any prime $p>2$ and a $p$-group $G$, then there exists an $(b(G) + 2)\times (b(G) + 2)$ upper triangular matrix $(A_{ij})_{j\geq i}$ such that:

\begin{enumerate}
    \item $A_{ij}\subseteq G$ are subgroups of $G$,
    \item $A_{11} = G$,
    \item $A_{ij} = \gamma_i(A_{1(j - i + 1)})$, $\forall j\geq i$,
    \item $A_{ij}\subseteq A_{i(j - 1)}$, $\forall j - 1\geq i$,
    \item $A_{ij}\supseteq A_{(i + 1)j}$, $\forall j\geq i + 1$,
    \item $A_{1(b(G) + 2)} = \{e\}$.
\end{enumerate}

\end{conjecture}

\begin{proposition}

The class-breadth conjecture for $p > 2$ is a consequence of Conjecture 1.

\end{proposition}

\begin{proof}
    From properties 2, 3, 5, 6 it follows that $\gamma_{b(G) + 2}(G) = \gamma_{b(G) + 2}(A_{11}) = A_{(b(G) + 2)(b(G) + 2)}\subseteq A_{1(b(G) + 2)} = \{e\}$.
\end{proof}

\begin{theorem}
    Consider any prime number $p > 2$. Conjecture 1 holds for all $p$-groups $G$ that satisfy:
    \begin{enumerate}
        \item[A1.] The group $H/[H, G]$ is non-cyclic for each normal non-abelian subgroup $H$ of $G$,

\item[A2.] $[G, H_2]\subseteq H_1$ for each normal subgroups $H_1, H_2$ of $G$ such that at least one of the following conditions holds:
    
    \begin{itemize}
        \item $H_1 = N_1'$, $H_2 = N_2'$ for some normal subgroups $N_1, N_2$ of $G$ with $N_1\subseteq N_2$ and $|N_2 : N_1| = p$,
        \item $H_2 = \gamma_i(H)$, $H_1 = \gamma_{i + 1}(H)$ for some normal subgroup $H$ of $G$ and some integer $i\geq 2$.
    \end{itemize}
    \end{enumerate}
\end{theorem}

\begin{proof}

We construct the $(b(G) + 2)\times (b(G) + 2)$ upper triangular matrix that satisfies properties 1 -- 6 of Conjecture 1 inductively as follows:

\begin{itemize}
    \item Define the $1\times 1$ matrix $(A_{11}) := (G)$ and let $C_{11} := C_{12}$ -- is any maximal subgroup of $G$.
    \item Assume that the $k\times k$ upper triangular matrix $(A_{ij})_{j\geq i}$ is defined for some integer $1\leq k < b(G) + 2$ and satisfies properties 1 -- 5 of Conjecture 1, and, moreover, satisfies:
    
    \begin{enumerate}
        \item[P1.] All subgroups $A_{ij}, j\geq i$ are normal in $G$.
        \item[P2.] If $k > 1$, then:
        
        \begin{enumerate}
            \item $A_{1j}/A_{1(j + 1)}$ is an abelian group, $\forall 1\leq j\leq k - 1$,
            \item $A_{1j}'/A_{1(j + 1)}'$ is an abelian group with $\log_p|A_{1j}'/A_{1(j + 1)}'|\leq b(G) - j + 1$, $\forall 1\leq j\leq k - 1$.
        \end{enumerate}
        \item[P3.] For each integer $1\leq j\leq k$ such that $A_{1j}$ is non-abelian, there exist normal subgroups $C_{j1}, C_{j2}$ of $G$ such that $C_{j1}\cup C_{j2}\subsetneq A_{1j}$, $|A_{1j} : C_{j1}| = |A_{1j} : C_{j2}| = p$, and $b_{A_{1j}}(x)\leq b(G) - j + 1$, $\forall x\in A_{1j}\setminus (C_{j1}\cup C_{j2})$.

        \item[P4.] If $k > 1$, then for each integer $1\leq j\leq k - 1$ such that $A_{1j}$ is non-abelian, there exist normal subgroups $N_{j1}, N_{j2}$ of $G$ such that $A_{1j}' = N_{j1}'$, $A_{1(j + 1)}' = N_{j2}'$, $N_{j2}\subseteq N_{j1}$, and $|N_{j1} : N_{j2}| = p$.
    \end{enumerate}

    \item Define normal subgroups $A_{i(k+1)}$, $1\leq i\leq k + 1$ of $G$, so that the $(k + 1)\times (k + 1)$ upper triangular matrix $(A_{ij})_{j\geq i}$ satisfies properties 1 -- 5 of Conjecture 1, satisfies properties P1 -- P4 (with $k + 1$ instead of $k$), and, in the case $k = b(G) + 1$, also satisfies property 6 of Conjecture 1, as follows:
    
    Define $H$ as the minimal under the inclusion normal subgroup of $G$ such that $H\subseteq A_{1k}$, $H\nsubseteq C_{k1}\cup C_{k2}$, and $H' = A_{1k}'$ (where $C_{j1}, C_{j2}, 1\leq j\leq k$ are defined as in property P3 for $k\times k$ submatrix of $(A_{ij})_{j\geq i}$).
    
    Consider the case when $H$ is a non-abelian subgroup of $G$. Let by definition $A_{1(k + 1)}$ be any normal subgroup of $G$ such that $A_{1(k + 1)}\subseteq H$, $|H : A_{1(k + 1)}| = p$, $A_{1(k + 1)}\nsubseteq C_{k1}\cup C_{k2}$, and $A_{1(k + 1)}\cap C_{k1} = A_{1(k + 1)}\cap C_{k2}$ (such a group exists from property A1 of Theorem 1 and Lemma 4). Therefore, from minimality of $H$, $A_{1(k + 1)}'\not= A_{1k}'$. Let $P := \{g\in A_{1(k + 1)}\:\vert\: [g, H]\subseteq A_{1(k + 1)}'\}$. $P$ is a normal subgroup of $G$, since $P = A_{1(k + 1)}\cap\pi^{-1}(Z_{\pi(G)}(\pi(H)))$, where $\pi : G\to G/A_{1(k + 1)}'$ is the homomorphism of factorization. In the case when $P = A_{1(k + 1)}$, it follows that $[A_{1(k + 1)}, H]\subseteq A_{1(k + 1)}'$, so $A_{1(k + 1)}/A_{1(k + 1)}'$ is a central subgroup of index $p$ in $H/A_{1(k + 1)}'$, i.e. $H/A_{1(k + 1)}'$ is an abelian group and $A_{1k}' = H' = A_{1(k + 1)}'$, contradiction. So $P\not= A_{1(k + 1)}$. From Lemma 1 there exists a normal subgroup $C_{(k + 1)1}$ of $G$ such that $P\subseteq C_{(k + 1)1}\subsetneq A_{1(k + 1)}$ and $|A_{1(k + 1)} : C_{(k + 1)1}| = p$. From Lemma 1, the fact that $A_{1(k + 1)}\cap C_{k1} = A_{1(k + 1)}\cap C_{k2}$, and the fact that $A_{1(k + 1)}\nsubseteq C_{k1}\cup C_{k2}$ it follows that there exists a normal subgroup $C_{(k + 1)2}$ of $G$ such that $A_{1(k + 1)}\cap (C_{k1}\cup C_{k2}) = A_{1(k + 1)}\cap C_{k1}\subseteq C_{(k + 2)2}\subsetneq A_{1(k + 1)}$ and $|A_{1(k + 1)} : C_{(k + 2)2}| = p$. Define $A_{i(k + 1)} := \gamma_i(A_{1(k + 2 - i)})$, $\forall 2\leq i\leq k + 1$.

    Consider the case when $H$ is an abelian subgroup of $G$. Then $A_{1k}' = H' = \{e\}$. Define $A_{1(k + 1)} := \{e\}$, $A_{i(k + 1)} := \gamma_i(A_{1(k + 2 - i)})$, $\forall 2\leq i\leq k + 1$.

    Proof of the fact that the $(k + 1)\times (k + 1)$ upper triangular matrix $(A_{ij})_{j\geq i}$ satisfies properties 1 -- 5 of Conjecture 1, satisfies properties P1 -- P4 (with $k + 1$ instead of $k$), and, in the case $k = b(G) + 1$, also satisfies property 6 of Conjecture 1:

    \begin{itemize}
        \item Properties 1 -- 3 are satisfied from the definition and the facts that $A_{1(k + 1)} = \gamma_1(A_{1(k + 1)})$, $A_{i(k + 1)} := \gamma_i(A_{1(k + 2 - i)})$, $\forall 2\leq i\leq k + 1$.

        \item Property 4 for $i = 1, j\leq k$ follows from property 4 applied to $k\times k$ submatrix of $(A_{ij})_{j\geq i}$.
        
        \item Property 4 for $i = 1, j = k + 1$ follows from the definition of $A_{1(k + 1)}$.

        \item Property 4 for $i\geq 2$ follows from property 4 for $i = 1$ and property 3, since $A_{ij} = \gamma_i(A_{1(j - i + 1)})\subseteq \gamma_i(A_{1(j - i)}) = A_{i(j - 1)}$.

        \item Property 5 for $j\leq k$ follows from property 5 applied to $k\times k$ submatrix of $(A_{ij})_{j\geq i}$.

        \item In the case $i = 1, j = k + 1$ and $H$ is a non-abelian group, property 5 is satisfied, because: from $A_{1(k + 1)}\subseteq H$, $H/A_{1(k + 1)}\cong\mathbb{Z}_p$, $H' = A_{1k}'$, and property 3 it follows that $A_{2(k + 1)} = A_{1k}' = H'\subseteq A_{1(k + 1)}$.

        \item In the case $i = 1, j = k + 1$ and $H$ is an abelian group, property 5 is satisfied, because: from $H' = A_{1k}'$ and property 3 it follows that $A_{2(k + 1)} = A_{1k}' = H' = \{e\}\subseteq A_{1(k + 1)}$.

        \item Property 5 holds for $i = 2$, $j = k + 1$, because: $k + 1 = j\geq i + 1 = 3$, so $k\geq 2$. From property P4 applied to $k\times k$ submatrix of $(A_{ij})_{j\geq i}$ it follows that there exist normal subgroups $N_{(k - 1)1}, N_{(k - 1)2}$ of $G$ such that $A_{1(k - 1)}' = N_{(k - 1)1}'$, $A_{1k}' = N_{(k - 1)2}'$, $N_{(k - 1)2}\subseteq N_{(k - 1)1}$, and $|N_{(k - 1)1} : N_{(k - 1)2}| = p$. Therefore, from the property A2 of Theorem 1 and property 3 it follows that $[A_{1(k - 1)}', G]\subseteq A_{1k}'$ and $A_{3(k + 1)} = [A_{2k}, A_{1(k - 1)}]\subseteq [A_{2k}, G] = [A_{1(k - 1)}', G]\subseteq A_{1k}' = A_{2(k + 1)}$.

        \item Property 5 holds for $i\geq 3$, $j = k + 1$, because: $k + 1 = j\geq i + 1$, so $k\geq i$. From property 3 it follows that $A_{(i - 1)k} = \gamma_{i - 1}(A_{1(k - i + 2)})$, $A_{i(k + 1)} = \gamma_i(A_{1(k - i + 2)})$, therefore, from $i\geq 3$ and property A2 of Theorem 1 it follows that $[A_{(i - 1)k}, G]\subseteq A_{i(k + 1)}$. From $k\geq i$, $[A_{(i - 1)k}, G]\subseteq A_{i(k + 1)}$, property 3, and property 5 applied to $k\times k$ submatrix of $(A_{ij})_{j\geq i}$ it follows that $A_{(i + 1)(k + 1)} = [A_{ik}, A_{1(k - i + 1)}]\subseteq [A_{ik}, G]\subseteq [A_{(i - 1)k}, G] \subseteq A_{i(k + 1)}$.

        \item Property 6 holds in the case $k = b(G) + 1$, because: from property P3 applied to $k\times k$ submatrix of $(A_{ij})_{j\geq i}$ and Lemma 5 it follows that $A_{1k} = A_{1(b(G) + 1)}$ is an abelian group. So $H\subseteq A_{1k}$ is also an abelian group and, from the definition of $A_{1(k + 1)}$, it follows that $A_{1(b(G) + 2)} = A_{1(k + 1)} = \{e\}$.

        \item Property P1 follows from the definition of groups $A_{ij}$.

        \item Property P2(a) is satisfied, because: from properties 3, 5 it follows that $A_{1j}' = A_{2(j + 1)}\subseteq A_{1(j + 1)}$, $\forall 1\leq j\leq k$.

        \item Property P2(b) for $1\leq j\leq k - 1$ follows from property P2(b) applied to $k\times k$ submatrix of $(A_{ij})_{j\geq i}$.

        \item In the case $j = k$ and $H$ is an abelian group, property P2(b) is satisfied, because: $A_{1k}' = H' = \{e\}$, so $A_{1k}'/A_{1(k + 1)}' = \{e\}$ and $\log_p|A_{1k}'/A_{1(k + 1)}'| = 0\leq b(G) - j + 1 = b(G) - k + 1$.
        
        \item In the case $j = k$ and $H$ is a non-abelian group, property P2(b) is satisfied, because: consider any element $g\in H\setminus (A_{1(k + 1)}\cup C_{k1}\cup C_{k2})$ (such an element exists, since $p > 2$, $H\nsubseteq C_{k1}\cup C_{k2}$ and, therefore, $H$ cannot be covered by its proper subgroups $A_{1(k + 1)}, C_{k1}\cap H, C_{k2}\cap H$). From $H' = A_{1k}'$ and \cite[Lemma 2.2]{10}, applied to the group $H$, its maximal subgroup $A_{1(k + 1)}$, and the element $g\in H\setminus A_{1(k + 1)}$, it follows that $\log_p|A_{1k}'/A_{1(k + 1)}'| = \log_p|H'/A_{1(k + 1)}'|\leq b_H(g)$. Therefore, from $k < b(G) + 2$, $A_{1(k + 1)}\subseteq H\subseteq A_{1k}$, $g\in H\setminus (A_{1(k + 1)}\cup C_{k1}\cup C_{k2})$, and property P3 applied to $k\times k$ submatrix of $(A_{ij})_{j\geq i}$ it follows that $\log_p |A_{1k}'/A_{1(k + 1)}'|\leq b_H(g)\leq b_{A_{1k}}(g)\leq b(G) - k + 1$. From $H/A_{1(k + 1)} \cong \mathbb{Z}_p$ it follows that $A_{1k}' = H'\subseteq A_{1(k + 1)}$, i.e. $A_{1k}'/A_{1(k + 1)}'$ is an abelian group with $\log_p |A_{1k}'/A_{1(k + 1)}'|\leq b(G) - k + 1$.

        \item In the case $j\leq k$, property P3 follows from property P3 applied to $k\times k$ submatrix of $(A_{ij})_{j\geq i}$.

        \item In the case $j = k + 1$ and $H$ is an abelian group, property P3 is satisfied, because: from the definition of $A_{1(k + 1)}$, $A_{1j} = A_{1(k + 1)} = \{e\}$ is an abelian group.

        \item In the case $j = k + 1$ and $H$ is a non-abelian group, property P3 is satisfied, because: $C_{(k + 1)1}\cup C_{(k + 1)2}\subsetneq A_{1(k + 1)}$ ($A_{1(k + 1)}$ cannot be covered by its two proper subgroups $C_{(k + 1)1}, C_{(k + 1)2}$). From the definition of $C_{(k + 1)1}, C_{(k + 1)2}, A_{1(k + 1)}, P$ it follows that $A_{1(k + 1)}\cap (P\cup C_{k1}\cup C_{k2})\subseteq C_{(k + 1)}\cup C_{(k + 1)2}\subsetneq A_{1(k + 1)}$. From $A_{1(k + 1)}\subseteq H\subseteq A_{1k}$ it follows that $$\{g\in A_{1(k + 1)}\:\vert\: b_{A_{1(k + 1)}}(g) \geq b_{A_{1k}}(g)\} = \{g\in A_{1(k + 1)}\:\vert\: b_{A_{1(k + 1)}}(g) = b_{A_{1k}}(g)\} = $$ $$ = \{g\in A_{1(k + 1)}\:\vert\: [g, A_{1(k + 1)}] = [g, A_{1k}]\} \subseteq P,$$ i.e.
            
            \begin{equation*}
                A_{1(k + 1)}\cap (\{g\in A_{1(k + 1)}\:\vert\: b_{A_{1(k + 1)}}(g) \geq b_{A_{1k}}(g)\}\cup C_{k1}\cup C_{k2})\subseteq
            \end{equation*}
            \begin{equation}
                    \subseteq A_{1(k + 1)}\cap (P\cup C_{k1}\cup C_{k2})\subseteq C_{(k + 1)}\cup C_{(k + 1)2}\subsetneq A_{1(k + 1)}.\tag{1.1}
            \end{equation}
            
            Therefore, from $(1.1)$, the fact that $A_{1(k + 1)}\subseteq A_{1k}$, and property P3 applied to $k\times k$ submatrix of $(A_{ij})_{j\geq i}$ it follows that $b_{A_{1(k + 1)}}(x)\leq b_{A_{1k}}(x) - 1\leq b(G) - k$, $\forall x\in A_{1(k + 1)}\setminus (C_{(k + 1)1}\cup C_{(k + 1)2})$.

            \item In the case $j\leq k - 1$, property P4 follows from property P4 applied to $k\times k$ submatrix of $(A_{ij})_{j\geq i}$.

            \item In the case $j = k$ and $H$ is a non-abelian group, property P4 is satisfied, because: from $H' = A_{1k}'$ and $|H : A_{1(k + 1)}| = p$ it follows that $N_{k1} := H$, $N_{k2} := A_{1(k + 1)}$ are suitable normal subgroups of $G$.

            \item In the case $j = k$ and $H$ is an abelian group, property P4 is satisfied, because: from $A_{1k}' = H' = \{e\}$, $A_{1k}$ is an abelian group.
    \end{itemize}
\end{itemize}
\end{proof}

\begin{remark}
    From \cite[Theorem 1.2]{10} and the fact that $A'\subseteq B$ for each $p$-groups $B\subseteq A$ with $|A : B| = p$ it follows that $H_1\subseteq H_2\subseteq G'$ and $H_2/H_1$ is an abelian group with $\log_p|H_2/H_1|\leq b(G)(b(G) + 1)/2$ for each normal subgroups $H_1, H_2$ of $G$ that satisfy property A2 of Theorem 1.
\end{remark}


\begin{thebibliography}{13}

\bibitem{6}
\author{\sc N. Blackburn}.
\textit{\em On a special class of p-groups}. Acta Math. 100 (1958).

\bibitem{3}
\author{\sc B. Eick, M.F. Newman, E.A. O'Brien}.
\textit{\em The class-breadth conjecture revisited}, Journal of Algebra, 300, 1, (2006).

\bibitem{2}
\author{\sc W. Felsch, J. Neubüser, W. Plesken}.
\textit{\em Space groups and groups of prime-power order, IV. Counterexamples to the class-breadth conjecture}, J. London Math. Soc. 24, 2, (1981).

\bibitem{4}
\author{\sc J. A. Gallian}.
\textit{\em On the breadth of a finitep-group}, Math Z 126, 224–226 (1972).

\bibitem{7}
\author{\sc H. G. Knoche}.
\textit{\em Über den Frobenius' schen Klassenbegriff in nilpotenten Gruppen}. Mathematische Zeitschrift 55 (1951).

\bibitem{1}
\author{\sc C. R. Leedham-Green, P. M. Neumann, J. Wiegold}.
\textit{\em The Breadth and the Class of a Finite P-Group}. Journal of the London Mathematical Society, s2-1: 409-420 (1969).

\bibitem{8}
\author{\sc I. D. Macdonald}.
\textit{\em Groups of breadth four have class five}. Glasgow Mathematical Journal 19, no. 2 (1978).

\bibitem{9}
\author{\sc V. D. Mazurov, E. I. Khurkhro (Editors)}.
\textit{\em Unsolved problems in group theory - The Kourovka notebook, No. 18}. arXiv:1401.0300v10, 8 Sep 2017.

\bibitem{5}
\author{\sc P. M. Neumann}.
\textit{\em An improved bound for BFCp-groups}. Journal of the Australian Mathematical Society 11, no. 1 (1970).

\bibitem{13}
\author{\sc P. Shumyatsky}.
\textit{\em Applications of Lie ring methods to group theory}. Nonassociative algebra and its applications. CRC Press., 373-396 (2019).

\bibitem{10}
\author{\sc A. A. Skutin}.
\textit{\em Proof of a conjecture of Wiegold}. Journal of Algebra 526 (2019).

\end{thebibliography}
\end{document}